\documentclass[10pt, conference]{IEEEtran}
\usepackage{graphicx}
\usepackage{amsmath}
\usepackage{cite}
\usepackage{url}
\usepackage{amsfonts,amssymb}
\usepackage{amsthm}
\IEEEoverridecommandlockouts

\begin{document}

\title{Smart Charging and Parking of Plug-in Hybrid Electric Vehicles in Microgrids Considering Renewable Energy Sources}
\author{Zheming~Liang,~\IEEEmembership{Student Member,~IEEE}
        and~Yuanxiong~Guo,~\IEEEmembership{Member,~IEEE}
\thanks{The authors are with the School
of Electrical and Computer Engineering, Oklahoma State University, Stillwater,
OK 74078 USA (e-mail: zheming.liang@okstate.edu; richard.guo@okstate.edu).}
}
\maketitle
\begin{abstract}
With the current trend of transforming a centralized power system into a decentralized one for efficiency, reliability, and environment reasons, the concept of microgrid that integrates a variety of distributed energy resources (DERs) on the distribution network is gaining popularity. In this paper, we investigate the smart charging and parking of plug-in hybrid electric vehicles (PHEVs) in microgrids with renewable energy sources (RES), such as solar panels, in grid-connected mode. To address the uncertainties associated with RES power output and PHEVs charging condition in the microgrid, we propose a two-stage scenario-based stochastic optimization approach with the objective of providing a proper scheduling for parking and charging of PHEVs that minimizes the average total operating cost while maintaining the reliability of the microgrid. A case study is conducted to show the effectiveness of the proposed approach. Extensive simulation results show that the microgrid can minimize the operating cost and ensure its reliability.
\end{abstract}

%----------------------------------------------------------------------------------
\section{Introduction}\label{sec:intro}
%----------------------------------------------------------------------------------

Driven by the goal of reducing greenhouse gas emissions, plug-in hybrid electric vehicles (PHEVs) have attracted lots of interest in recent years. However, recent research has shown that even a modest penetration of PHEVs could cause frequent overloads in the distribution system if they are allowed to charge whenever they are plugged into the system~\cite{SaOr15}, since this action may increase the peak loads. In order to avoid such problems, renewable energy sources (RES) such as solar panels are needed. The solar panels can support the peak demand in the peak periods, while combined heat and power (CHP) units can support power and heat demand at off-peak periods. However, bulk power systems are not designed to manage large numbers of small distributed energy sources (DERs) such as CHP and RES, therefore, the concept of microgrid has been proposed~\cite{HaAs07}. According to the U.S. Department of Energy~\cite{Mi15}, a microgrid is ``a group of interconnected loads and DERs within clearly defined electrical boundaries that acts as a single controllable entity with respect to the grid and that connects and disconnects from such grid to enable it to operate in both grid-connected or isolated mode''. The microgrid enables a more flexible bulk power system by handling deferrable loads and the RES locally, which can be used to resolve the overload issues.

In summary, the main contributions of this paper are as follows.
\begin{enumerate}
    \item We propose a two-stage scenario-based stochastic optimization formulation to model the energy management problem in microgrids with CHP and RES which minimizes the expected operating cost of the grid-connected microgrid considering the charging and parking problems of the PHEVs.

    \item We transform the two-stage stochastic optimization problem into a large-scale mixed integer linear program (MILP). The uncertainties related with the RES power output and PHEVs charging condition are modeled by scenarios.

    \item We conduct a case study based on real-world datasets to present the effectiveness of our proposed approach. Simulation results show the proposed model is reliable in various scenarios.
\end{enumerate}
The rest of the paper is organized as follows. Section~\ref{sec:model} describes the system model and the assumptions we use in this paper. Section~\ref{sec:prob} presents the mathematical formulation for the two-stage stochastic optimization problem. In Section~\ref{sec:solu}, we propose an efficient solution method to solve the optimization problem. We conduct a case study in Section~\ref{sec:case} and then draw conclusions in Section~\ref{sec:con}.
\pagenumbering{arabic}
\setcounter{page}{1}
%----------------------------------------------------------------------------------
\section{System Modeling}\label{sec:model}
%----------------------------------------------------------------------------------

We consider a microgrid consisting of combined heat and power (CHP), renewable energy sources (RES), PHEVs, deferrable loads, and base power and heat loads. The energy management is over a finite time horizon $\mathcal{T} := \{1, 2, \ldots, T\}$ (e.g., 24 hours) with time period indexed by $t$.

%-------------------------------------
\subsection{Combined Heat and Power}
%-------------------------------------

The power output of each CHP unit are restricted by its generator capacity:
\begin{equation}\label{con:chp_p_limit}
\underline{P}_{i} \leq p_{i, t}^s \leq \overline{P}_{i}, \forall i, t, s
\end{equation}
where $\underline{P}_{i}$ and $\overline{P}_{i}$ are the minimum and maximum power output for the $i$-th CHP unit.

%-------------------------------------
\subsection{Plug-in Hybrid Electric Vehicles}
%-------------------------------------

We use $E_{m, t}^s$ to denote the energy storing in $m$-th PHEV at the end of the period $t$, with initial available energy $E_{m, 0}^s$. Then, we have the following dynamics for the stored energy in the $m$-th PHEV:
\begin{equation}\label{con:es_dynamics}
E_{m, t}^s = E_{m, t - 1}^s + r_{m, t}^{s, +} \eta_{m}^{+} - r_{m, t}^{s, -}/\eta_{m}^{-}, \forall m, t, s
\end{equation}
where $r_{m, t}^{s, +}$ and $r_{m, t}^{s, -}$ are the power charged into or discharged from the $m$-th PHEV at time $t$ in scenario $s$, and $\eta_{m}^{+}$ and $\eta_{m}^{-}$ represent the charging and discharging efficiencies of the $m$-th PHEV, respectively.

Each PHEV has a finite capacity, and therefore the stored energy in it must have the following lower and upper bounds:
\begin{equation}
\underline{E}_m \leq E_{m, t}^s \leq \overline{E}_{m}, \forall m, t, s
\end{equation}
where the upper bound $\overline{E}_{m}$ is the storage capacity for $m$-th PHEV and the lower bound $\underline{E}_m$ is imposed to reduce the impact of deep discharging on the storage lifetime (e.g., 20\% of the capacity for Li-ion batteries).

Furthermore, PHEVs have charging or discharging rate limit as follows:
\begin{equation}
0 \leq r_{m, t}^{s, +} \leq R_{m}^{+} I_{m,t}^s, 0 \leq r_{m, t}^{s, -} \leq R_{m}^{-} I_{m,t}^s, \forall m, t, s
\end{equation}
where $R_m^{+}$ and $R_m^{-}$ denote the maximum charged and discharged energy over period $t$ for the $m$-th PHEV, respectively. We use a binary variable $I_{m,t}^s$ to represent the parking condition. When $I_{m,t}^s$ is 1, then the $m$-th PHEV is parked at the charging station, which is available to charge or discharge. When $I_{m,t}^s$ is 0, the $m$-th PHEV is on the road or not able to charging and discharging.

Lastly, without loss of generality, we assume that for each PHEV, the final stored energy should be the same as the initial stored energy, i.e.,
\begin{equation}
E_{m, T}^s = E_{m, 0}^s, \forall m, s.
\end{equation}

%-------------------------------------
\subsection{Energy Loads}
%-------------------------------------

The microgrid contains both power and heat energy demand. In this paper, power loads are classified into two categories: static power loads and deferrable power loads. We denote the aggregate power demand of all static loads as $P_t^0$ that must be satisfied at each period $t$.

For deferrable power loads, they only require a certain amount of electric energy to be delivered over a specified time interval and therefore, have some flexibility in their power profiles, such as electric water heaters. The loads are required to be served within the time interval $[T_j^{a}, T_j^{d}]$. The energy $L_j^s$ must be served with a minimum and a maximum load serving rate of $\underline{l}_j$ and $\overline{l}_j$, respectively. The load serving requirement is expressed as:
\begin{equation}\label{con:load}
\sum_{t = T_j^{a}}^{T_j^{d}} l_{j, t}^s = L_j^s, \underline{l}_j \leq l_{j, t}^s \leq \overline{l}_j, \forall j, t \in [T_j^a, T_j^d], s
\end{equation}
\begin{equation}\label{con:load2}
l_{j, t}^s = 0, \forall j, t \notin [T_j^a, T_j^d], s
\end{equation}
where $l_{j, t}^s$ is the power delivered to the $j$-th task over period $t$ for scenario $s$.

%----------------------------------------------------
\subsection{Power Exchange with the Main Grid}
%----------------------------------------------------

In grid-connected mode, the microgrid can either import or export power from or into the main grid. As with~\cite{ZhGa13}, we assume that for the microgrid in the grid-connected mode, the surplus power can be sold to the main grid with known selling price $c_t^{+}$, and the shortage power can be purchased from the main grid with known purchase price $c_t^{-}$ at each time period $t$. We further denote $g_t^{s, +}$ and $g_t^{s, -}$ as the amount of power sold into and bought from the main grid at time $t$, respectively. Then, the market exchange cost of the microgrid with the main grid at time $t$ can be represented as $c_t^{+} g_t^{s, +} - c_t^{-} g_t^{s, -}$. Moreover, we use $\overline{g}_t$ to represent the capacity limit in~\eqref{eq:connection1}. Then, we have:
\begin{equation}\label{eq:connection1}
0 \leq g_t^{s, -} \leq \overline{g}_t, 0 \leq g_t^{s, +} \leq \overline{g}_t, \forall t, s.
\end{equation}

%----------------------------------------------------
\subsection{Energy Balance}
%----------------------------------------------------

Both power and heat must be balanced in all time periods for the microgrid. Since the CHP can generate power and useful heat simultaneously, we use a power-to-heat ratio $\alpha_i$ to represent the relationship between power and useful heat output of CHP unit $i$. That is, when one unit of power is generated, CHP unit $i$ generates $\alpha_i$ unit of useful heat. Also, $N_c$, $N_{p}$, and $N_j$ are adopted to represent the CHP units, PHEVs, and deferrable loads, respectively. Therefore, we have the following power and heat balance equations:
\begin{align}\label{con:p_balance}
& \sum_{i = 1}^{N_c} p_{i, t}^s + \sum_{m = 1}^{N_{p}} \left[r^{s, -}_{m, t}- r^{s, +}_{m, t}\right] + w_t^s  \nonumber \\
& = g_t^{s, -} - g_t^{s, +} + P_t^{0} + \sum_{j = 1}^{N_j} l_{j, t}^s, \forall t, s.
\end{align}
Note that in the above equation, $w_t$ denotes the aggregate RES power output in the microgrid and therefore, is uncertain.
\begin{equation}\label{con:h_balance}
\sum_{i = 1}^{N_c} \alpha_i p_{i,t}^s \geq Q_t, \forall t, s
\end{equation}
where we assume that all heat demand must be satisfied locally due to the fact that heat cannot be transported over a long distance, and surplus heat can be disposed without penalty.

%-------------------------------------------
\section{Two-Stage Stochastic Formulation}\label{sec:prob}
%-------------------------------------------

The objective of the proposed framework is to minimize the expected operating cost of the microgrid , which is shown as follows:
\begin{align}\label{eq:twostage}
& \min \sum_{i = 1}^{N_{c}} \sum_{t=1}^{N_t} \sum_{s=1}^{N_s} \rho_s c_i^p p_{i, t}^s \nonumber \\
& \quad \quad + \sum_{m = 1}^{N_{p}} \sum_{t=1}^{N_t} \sum_{s=1}^{N_s} \rho_s c_{m}^p \big(r_{m, t}^{s, +} \eta_m^{+} + r_{m, t}^{s, -}/\eta_m^{-} \big) \nonumber \\
& \quad \quad - \sum_{t=1}^{N_t} \sum_{s=1}^{N_s} \rho_s c_t^{-} g_t^{s, -} + \sum_{t=1}^{N_t} \sum_{s=1}^{N_s} \rho_s c_t^{+} g_t^{s, +},
\end{align}
subject to constraints~\eqref{con:chp_p_limit}--constraints~\eqref{con:h_balance}. $N_t$ and $N_s$ represent the total time horizon and the total scenarios. $c_m^p$ is the degradation cost of the $m$-th PHEV due to charging and discharging, respectively. Note that the first line represents the operation and maintenance cost of CHP units, and the second line shows the operating cost of PHEVs, and the last line denotes the cost or revenue when exchanging power with the main grid.

%-------------------------------------------
\section{Solution Approach}\label{sec:solu}
%-------------------------------------------

We use a scenario-based two-stage stochastic programming to handle the uncertainties relate with the RES power output and the PHEVs charging condition. The uncertainties are modeled by scenarios. We assume that the historic data for the RES power output can be obtained by the microgrid central controller and the PHEVs charging condition follows Bernoulli distribution. Also, the uncertainties are assumed to be independent from each other. Then we use the historic data to generate 3000 scenarios with even probability. Moreover, we use a fast-forward scenario reduction method~\cite{DuJi03} to reduce the original 3000 scenarios into 25 scenarios. The objective of the scenario reduction is to obtain a small set of scenarios that can maintain properties of original scenarios at a reasonable level. Therefore, the microgrid energy management problem is a large-scale mixed integer linear program (MILP), which can be directly solved by commercial solvers such as GUROBI 6.0.0~\cite{Gu15}.

%-------------------------------------------
\section{Case Study}\label{sec:case}
%-------------------------------------------

In this section, we present the results of a case study by evaluating the proposed algorithm using real-world datasets. Here, we first describe the details of the data we used for modeling the microgrid. Then, we evaluate the performance of our stochastic solution.

All simulations are implemented on a desktop computer with 4.0 GHz Intel Core i7-4790 CPU and 8GB RAM. The scenario generation and scenario reduction are done by MATLAB~\cite{Ma15}. The microgrid energy management problem are simulated with Python~\cite{Python}.

%%%%%%%%%%%%%%%%%%%%%%%%%%%%%%
\subsection{Microgrid Description}
%%%%%%%%%%%%%%%%%%%%%%%%%%%%%%

\begin{table}[!t]
\centering
\caption{PHEV parameters}
\label{table:3}
\begin{tabular}{c c c c}
\cline{1-4}
$\underline{E}_m$ (kWh) & $\overline{E}_{m}$ (kWh) & $R_{m}^{p, +}$ (kW) & $R_{m}^{p, -}$ (kW) \\
\hline
4 & 18 & 4 & 4    \\
\hline
\end{tabular}
\end{table}
The considered microgrid includes three CHP units, several solar panels, several PHEVs, several electric water heaters, a base power load, and a base heat load. The parameters of these three CHP units are listed in~\cite{ChGo12}. The heat to power ratio is set to 1.2. The total installed capacity of the solar panels is 1000 kW. The historic data for solar panel power output $w_t$ is taken from~\cite{CoCl12} where coefficients are appropriately scaled. As shown in Table~\ref{table:3}, we have 50 identical Chevrolet 2016 PHEVs~\cite{Volt2016}. Each PHEV battery has a storage capacity of 18 kWh, and has charging and discharging efficiencies of 0.9. The maximum charging and discharging rates for each battery are both set to be 4 kW. To prolong the battery lifetime, the energy level should not drop below $20\%$ of its capacity. Both the initial and final stored electrical energy are set to be $50\%$ of the total capacity. The battery degradation cost coefficient $c_m^{\text{es}}$ is set to be 0.0035 \$/kWh using the calculation method in~\cite{KeTo05B}. We use the historic dataset for the forecasted total power demand~\cite{WaSh08} where coefficients are appropriately scaled. The electricity prices are taken from real-world electricity prices~\cite{PJM15} where coefficients are appropriately scaled. The selling price is set to be $80\%$ of the purchasing price. The electricity transmission limit between the microgrid and the main grid is set to be $4000 kW$. The data for base power and heat loads are taken from~\cite{Solar}. The data for electric water heaters are taken from~\cite{SaOr15}.

%%%%%%%%%%%%%%%%%%%%%%%%%%%%%%%%%%%%%%%%
\subsection{Results and Discussion}
%%%%%%%%%%%%%%%%%%%%%%%%%%%%%%%%%%%%%%%%

\subsubsection{Stochastic Approach versus Deterministic Approach}

\begin{figure}[!t]
\centering
\includegraphics[width=2.5in]{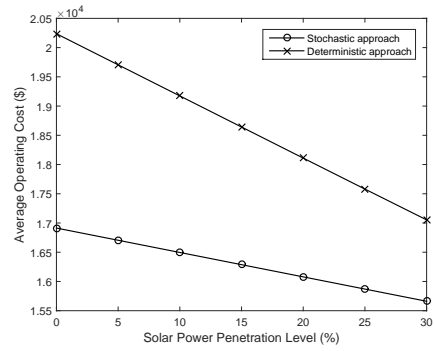}
\caption{The average operating cost with the increasing solar power penetration level in stochastic approach and in deterministic approach.}
\label{fig:Solar_change}
\end{figure}
In this part, we test our proposed framework in two approaches, stochastic approach and deterministic approach. As shown in Fig.~\ref{fig:Solar_change}, with the increase of the penetration level of solar power in the microgrid, the average operating cost of both approaches are decreasing. It is because that the solar panels can provide power to mitigate the power demand from the PHEVs and the base load in the peak periods with no operating cost. Moreover, we also observed that the stochastic approach has a much lower average operating cost than that of the deterministic approach. It is because that the stochastic approach is based on the scenarios that are capturing the characteristics of all original 3600 scenarios. However, the deterministic approach can only take one scenario which can not guarantee to have the optimal average operating cost. Therefore, the stochastic approach is more capable in handling uncertainties.

\subsubsection{Sensitivity Analysis}

\begin{figure}[!t]
\centering
\includegraphics[width=2.5in]{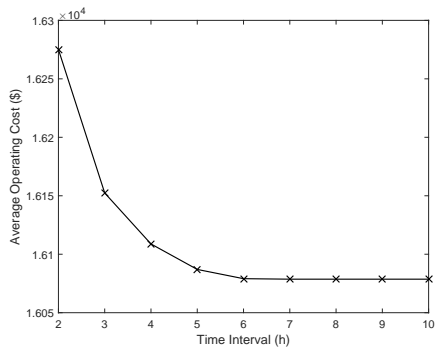}
\caption{The average operating cost with the increasing deferrable load serving time interval in stochastic approach.}
\label{fig:Time_change}
\end{figure}
In this part, we conduct a test our proposed model to observe the impact of the length of the time interval on the microgrid average operating cost. As shown in Fig.~\ref{fig:Time_change}, with the increase of the deferrable load serving time interval, the microgrid has more and more flexibility in system operating, therefore the microgrid average operating cost decreases fast at beginning and becomes stable at last. It is because the deferrable loads, such as electric water heaters can be served in the periods with lower electricity price when the time interval is large. However, when the time interval is small, the deferrable loads are forced to be served within such periods regardless of the electricity price. Therefore, the proper deferrable load serving time intervals are important for an optimal microgrid energy management problem.

%----------------------------------------
\section{Conclusions}\label{sec:con}
%----------------------------------------

We have proposed a framework for smart charging and parking of plug-in hybrid electric vehicles (PHEVs) in microgrids with renewable energy sources (RES) and combined heat and power (CHP) units. In our model, uncertainties are from the renewable energy sources power output and the PHEVs charging condition. We formulated a two-stage scenario based stochastic program to handle the uncertainties in the microgrid. Extensive simulations have been conducted to demonstrate the effectiveness of our proposed framework.

\bibliographystyle{IEEEtran}
\bibliography{./myref}

\end{document}